\documentclass[10 pt, letterpaper]{article}

\usepackage{amsmath}
\usepackage{graphicx}

\newcommand{\ph}{\phantom}

\begin{document}

\centerline{\textbf{Some Knots With Surgeries Yielding Lens Spaces}}
\medskip
\medskip
\centerline{John Berge}

\begin{abstract}
This is a facsimile of the circa 1990 unpublished manuscript with the same title. All the original text, figures and tables are included; although text has been reset in \TeX, the original hand-drawn figures have been redrawn digitally, and the parameter $k$ in the original table of lens spaces has been replaced with the originally intended $\lambda$. And, of course, this abstract has been added.
\end{abstract}

\begin{flushright}
Version 1.0
\end{flushright}

\medskip

We describe several families of knots in $S^3$ which have non-trivial Dehn-surgeries
yielding lens spaces. These knots are characterized by Theorem 1 and appear to
include all known knots in $S^3$ which have surgeries yielding lens spaces.

\underline{Theorem 1}.
Let $k$ be a simple closed curve on the boundary of a standardly embedded genus two
handlebody $H$ in $S^3$. If $k$ represents a free generator of both the fundamental group of $H$ and the fundamental group of the complementary handlebody $H^\prime$ in $S^3$, then there is an integer surgery on $k$ which yields a lens space.

\underline{Proof}.
Observe that because $k$ represents a free generator of $\pi_1(H)$, there exists a set
of cutting disks $D_1$, $D_2$ for $H$ with the property that $k$ intersects $D_1$ exactly once transversely, and $k$ is disjoint from $D_2$. Similarly, there exists a set of cutting disks $D_1^\prime$, $D_2^\prime$ for $H^\prime$ with the property that $k$ intersects $D_1^\prime$ exactly once transversely, and $k$ is disjoint from $D_2^\prime$. Note that adding the 2-handle $D_2^\prime$ to $H$ yields a manifold with torus boundary embedded in $S^3$. Thus this manifold is the exterior of a knot in $S^3$. Note also that adding the 2-handle $D_1^\prime$ to the boundary of this manifold yields $S^3$. Hence $D_1^\prime$ is the meridian of this knot. Next, note that because $k$ is transverse to $D_1^\prime$ at a single point, surgery on this knot along $k$ is an integer surgery. Observe also, since $k$ intersects $D_1$ only once, that adding a 2-handle to $H$ along $k$ yields a solid torus. Thus surgery on this knot along $k$ yields a manifold which has a Heegaard splitting of genus one and which must therefore be a lens space. Finally, note because $k$ is transverse to $D_1^\prime$ at a single point, $k$ is isotopic to the core of this knot, and so we may take $k$ itself as representing the knot on which the surgery is performed. Q.E.D.

\underline{Remark}.
We will sometimes informally refer to a curve $k$ satisfying Theorem 1 as a 
\underline{double-primitive} and to the knot $k$ as a \underline{double-primitive knot}. We also mention that since the surgery curve $k$ has the same knot type as the knot itself, we will often blur the distinction between the surgery curve and the knot.

Given a knot $k$ in $S^3$ which has an integer surgery that yields a lens space, it is
in general not easy to determine whether $k$ satisfies Theorem 1. Of course, if
$k$ satisfies Theorem 1, $k$ must be a tunnel-number-one knot, but in general a tunnel-number-one knot will have many inequivalent tunnels and finding all of these tunnels and checking whether any of them lead to a decomposition of the exterior of $k$ such that the surgery curve on $k$ satisfies Theorem 1 is not an easy matter. Hence we sought to find more explicit descriptions of the knots in $S^3$ satisfying Theorem 1 and to find explicit formulas giving the lens spaces that result from surgeries on these knots. These descriptions and formulas should be useful in determining whether a knot with surgery yielding a lens space satisfies Theorem 1 and if it does satisfy Theorem 1 they should be useful in precisely identifying the knot.

We also show that the knots satisfying Theorem 1 correspond to knots in lens spaces which are either 0 or 1-bridge knots in lens spaces and which have integer surgeries yielding $S^3$. This alternative characterization may also prove useful in identifying knots of this type, and we describe it below.

We must also include a warning that there is a chance that the list of explicit descriptions
of the knots satisfying Theorem 1, which we give, is not complete. We believe that there is fairly good evidence that it is complete but we do not have a proof of this fact. 

\thanks{Thanks to Professor Jim Cannon and my major professor D.R.McMillan Jr. for their patience and support.}

This leaves open the possibility that there are knots that satisfy Theorem 1 which are not on any of the lists we give.

At the end of the paper, we have included a table which gives those lens spaces $L(p,q)$ with $p \leq 500$ that are obtainable via surgeries on the knots in our list. However, to keep the size of the table small, we have not included those lens spaces which are obtainable only via surgeries on torus knots or cables about torus knots.

\medskip
\centerline{\underline{Introductory Definitions and Lemmas}.}
\medskip

A \underline{Heegaard splitting} $(M,F)$, of a closed, connected, orientable 3-manifold $M$ is a closed, connected, orientable surface $F$ embedded in $M$ which separates $M$ into two components whose closures are handlebodies $V$ and $W$. If $(M,F)$ is a Heegaard splitting, then a choice of complete sets, say $v$ and $w$, of meridional disks for $V$ and $W$ yields a \underline{Heegaard diagram} $(F,\partial v,\partial w)$, of $(M,F)$. We will also use the notation $(F,v,w)$ for a Heegaard splitting when $M$, $V$ and $W$ are understood and our interest is in $F$ and sets of cutting disks 
$v$ for $V$ and $w$ for $W$. We always assume, unless indicated otherwise, that the number of intersections of $\partial v$ with $\partial w$ is minimal up to isotopy in $F$.

It is known that genus one Heegaard splittings of lens spaces are unique up to isotopy [].
Thus if $L(p,q)$ is a lens space, then $H_1(L(p,q))$ has a pair of distinguished generators which are represented by the cores of the two solid tori of a genus one Heegaard splitting of $L(p,q)$. Here, and throughout the paper, we use integral homology. Let $k$ be a knot in $L(p,q)$ and $\overrightarrow{g}$ one
of the distinguished generators of $H_1(L(p,q))$. Note that if $[\overrightarrow{k}]=\lambda[\overrightarrow{g}]$ and $(\lambda,p) = 1$, then $\lambda$ is determined up to $\pm \lambda^{\pm 1} \mod  p$. Thus, $\pm \lambda^{\pm 1} \mod p$ is an invariant of $k$. We will have use for this invariant later.

Let $M$ be a lens space. A knot $k$ in $M$ is a \underline{0-bridge knot} if $k$ is isotopic to a knot $k^\prime$ lying in a torus which is a Heegard surface of a genus one Heegaard splitting of $M$. A knot $k$ in $M$ is a \underline{1-bridge knot} in $M$ if $k$ is not a 0-bridge knot in $M$ and $k$ is isotopic to a knot $k^\prime$ which intersects each of the solid tori of a genus one Heegaard splitting of $M$ in a single unknotted arc.

If $k$ is a knot in a 3-manifold $M$, the \underline{exterior} of $k$ is the complement of an open tubular neighborhood of $k$ in $M$. Note that if $k$ is a knot in manifold $M$ and integer surgery on $k$ yields a manifold $N$, then there is an associated knot $k'$ in $N$ whose exterior is homeomorphic to the exterior of $k$ and which has an integer surgery yielding $M$.

\underline{Lemma 1}.
Let $k$ be a knot in $S^3$ which satisfies Theorem 1, suppose that integer surgery on $k$ yields a lens space $M$, and let $k^\prime$ be the associated knot in $M$, then $k^\prime$ is 0-bridge or 1-bridge in $M$.

\underline{Proof}.
We use the notation used in the proof of Theorem 1. Observe that when a 2-handle is added to $H$ along the surgery curve $k$, $H$ becomes a solid torus $V$ and $\partial V$ is a genus one Heegaard splitting surface for $M$. Let $V^\prime$ be the other solid torus of this genus one Heegaard splitting. Let $A$ be an annular neighborhood of the surgery curve $k$ in $\partial H$. Note $\partial D_1^\prime$ represents the 
associated knot $k^\prime$ in $M$ and $\partial D_1^\prime \cap A$ is a subarc of $\partial D_1^\prime$ which is an essential arc $\alpha$ in $A$. Let $\beta$ be the subarc arc of $\partial D_1^\prime$ that does not lie in $\text{int}(A)$. Then $\beta$ is an unknotted arc which lies in $\partial V$, and so it is certainly unknotted in $V^\prime$. 

Next consider the subarc $\alpha$ of $\partial D_1^\prime$. Let $\gamma$ denote the arc $\partial D_1 \cap A$. Then $\gamma$ is an essential arc in $A$. And we may assume that $A$ has been chosen so that $\alpha$ and $\gamma$ are disjoint.
Then $\alpha$ is parallel into the boundary of $V$ via an isotopy which pushes the interior of $\alpha$ across a disk in $H$ isotopic to $D_1$, while keeping the endpoints of $\alpha$ fixed. In other words, $\alpha$ is unknotted in $V$. Thus the knot $k^\prime$ is the union of two arcs, $\alpha$ and $\beta$, intersecting only in their endpoints, with $\alpha$ unknotted in $V$ and $\beta$ unknotted in $V^\prime$. So $k^\prime$ is 0-bridge or 1-bridge in $M$. Q.E.D.

\underline{Lemma 2}.
Let $(F,\partial v,\partial w)$ be a genus two Heegaard diagram of $S^3$ and suppose that $k$ is a simple closed curve in $F$ such that $k$ intersects $\partial v_1$ transversely exactly once, $k$ intersects $\partial w_1$ transversely exactly once and $k$ is disjoint from $\partial v_2$ and $\partial w_2$. Suppose each curve $\partial v_1,\partial v_2,\partial w_1$, and $\partial w_2$ is oriented, then all intersections
of $\partial v_2$ with $\partial w_2$ have the same sign. If, in addition, there is no edge of the diagram with its endpoints on $\partial v_1$ or $\partial w_1$, then the curves of the diagram can be oriented so that all intersections of $\partial v_n$ with $\partial w_m$ have the same sign for $n,m \in \{1,2 \}$. 

\underline{Proof}
First note that because $k$ intersects the cutting disks of each handlebody only once, there are no edges of the diagram $(F,\partial v,\partial w)$ which have both endpoints on the same side of $\partial v_2$ or $\partial w_2$. Next, note that if there are edges of the diagram with both endpoints on the same side of $\partial v_1$ or $\partial w_1$,
then we can replace $\partial v_1$ or $\partial w_1$ with a new cutting disk and obtain a new diagram with fewer intersections and still maintain the property that $k$ intersects the cuttings disks in each handlebody only once. By [] if this diagram is nontrivial, then it has a wave $\omega$. This now forces the claimed result as can easily be seen by examining the possible diagrams that can occur. Q.E.D.
  
Lemma 2 allows us to show that a 1-bridge knot in a lens space $M$, which has a Dehn-surgery yielding $S^3$, is embedded in $M$ in a particularly simple and pleasant manner.

\underline{Theorem 2}.
Let $(F, v, w)$ be a genus one Heegaard splitting of a lens space $M$ and let $k$ be a 1-bridge knot in $M$ such that Dehn-surgery on $k$ yields $S^3$. Assume that $\partial v$ and $\partial w$ intersect minimally up to isotopy in $F$. Then $k$ is isotopic to a knot $k^\prime$ in $M$ such that $k^\prime$ is the union of two arcs $\alpha$ and $\beta$, with $\alpha$ properly embedded in $v$ and $\beta$ properly embedded in $w$.

\underline{Proof}.
Since $k$ is 1-bridge in $M$, we may assume that $k$ is the union of two arcs $\alpha$ and $\beta$ with $\alpha$ unknotted in $V$ and $\beta$ unknotted in $W$. It follows that we may also assume that $\beta$ lies in $F$ and that $\alpha$ is properly embedded in $V$ with $\alpha$ disjoint from $v$. Note that if $\beta$ is not disjoint from $w$, then there is an isotopy of $w$ which makes $w$ disjoint from $\beta$.

Let $N$ be a small regular neighborhood of $\alpha$ in $V$. Then $V$ minus the interior of $N$ is a handlebody $H$ of genus two. Let $A$ be the annulus in which $N$ meets $H$, and let $m$ be an essential simple closed curve in $A$. Note that $m$ represents a primitive element in $\pi_1(H)$ and there is a cutting disk $v^\prime$ for $H$, disjoint from $v$, such that $\partial v^\prime$ intersects $m$ transversely exactly once.

Now consider doing surgery on $k$ to get $S^3$. Since $k$ is not 0-bridge, it follows that the exterior of $k$ is not Seifert-fibered. Since $M$ is irreducible, and surgery on $k$ yields $S^3$, the exterior of $k$ must be irreducible. Then the Cyclic-Surgery-Theorem of [C,G,L,S] shows that the surgery which yields $S^3$ must be an integer surgery.

Next, note that the exterior of $k$ is homeomorphic to $H \cup w$. Let $c$ be a simple closed curve on $\partial H$ to which a 2-handle is added to perform the surgery. Then $c$ intersects $m$ transversely exactly once. Thus we have a genus two Heegaard splitting of $S^3$ with $m$ intersecting the cutting disks of each handlebody once. Lemma 2 then implies that if $\partial v$ and $\partial w$ are oriented, the intersections of $\partial v$ and $\partial w$ all have the same sign. This in turn implies that $k$
is isotopic to a knot of the claimed form in $M$. Q.E.D.

\underline{Remark}.
Observe that Theorem 2 shows that a 1-bridge knot, in a lens space $M$, which has a surgery yielding $S^3$, is isotopic to a knot which is simultaneously braided with respect to both of the solid tori of a genus one Heegaard splitting of $M$.
  
The next theorem, which follows from Theorem 2, shows that 1-bridge knots in a lens space, which have surgeries yielding $S^3$, are distinguished by the homology classes that they represent.

\underline{Theorem 3}.
Suppose that $k$ and $l$ are \emph{1}-bridge knots in a lens space $M$ such that surgery on both $k$ and $l$ yields $S^3$ and $\overrightarrow{k}$ is homologous to $\pm \overrightarrow{l}$ in $M$. Then $\overrightarrow{k}$ is isotopic to $\pm\overrightarrow{l}$ in $M$.

\underline{Proof}
The proof is easy and follows from the fact that $k$ and $l$ are isotopic to knots with the form shown to exist in Theorem 2. Q.E.D.

The following formula allows one to determine the lens space that results from performing surgery along a surgery curve that represents a double primitive knot in Theorem 1, provided one knows the homology class that the surgery curve represents in the boundary of the handlebody $H$.

Choose a standard basis for $H_1(\partial H)$ as in Fig. 1 with $\overrightarrow{a},\overrightarrow{b}$ bounding disks in $H$ 
and $\overrightarrow{A},\overrightarrow{B}$ bounding disks in $H^\prime$.

\centerline{Insert Fig. 1 here.}

\underline{Lemma 3}.
Suppose $k$ is a simple closed curve in $\partial H$ satisfying Theorem 1. If 
\[ [\overrightarrow{k}] = A[\overrightarrow{A}]+B[\overrightarrow{B}] +a[\overrightarrow{a}]+b[\overrightarrow{b}] \]
in $H_1(\partial H)$, then surgery along $k$, on the knot represented by $k$, yields a lens space
$L(p,q)$, with
\begin{align*}
p &= |Aa+Bb|, \\
a^2q &\equiv \pm B^{\pm 2} \mod p, \text{ and } \\
b^2q &\equiv \pm A^{\pm 2} \mod p.
\end{align*}

\underline{Proof}.
Let $h$ be an orientation preserving homeomorphism of $H^\prime$ such that $h(\overrightarrow{a})=\overrightarrow{k}$.
Then $h$ induces a sympletic transformation $\phi$ on $H_1(\partial H^\prime)$. Since intersection numbers are preserved, we must have:
\begin{align*}
\phi(\overrightarrow{A}) &= (\ph{-}s,-r,\ph{-}0,\ph{-}0), \\
\phi(\overrightarrow{B}) &= (-b,\ph{-}a,\ph{-}0,\ph{-}0), \\
\phi(\overrightarrow{a}) &= (\ph{a}A,\ph{b}B,\ph{-}a,\ph{-}b), \\
\phi(\overrightarrow{b}) &= (\,\ph{-}i,\ph{-}j,\ph{-}r,\ph{-}s)
\end{align*}
where $i,j,r,$ and $s$ are integers satisfying $as-br=1$ and $Ar+Bs=ai+bj$. Now we want to change the basis of $H_1(H)$ so that the surgery curve $\overrightarrow{k}$ represents $(1,0)$ and project everything onto the other generator $(0,1)$ of $H_1(H)$. This can be done by forming the dot-product of $\phi(\overrightarrow{A})$, $\phi(\overrightarrow{B})$,
$\phi(\overrightarrow{a})$, and $\phi(\overrightarrow{b})$ with $(-B,A,0,0)$ and projecting to get:
\begin{align*}
\phi^\prime(\overrightarrow{A}) &= -Ar-Bs, \\
\phi^\prime(\overrightarrow{B}) &= \ph{-}Aa+Bb, \\
\phi^\prime(\overrightarrow{a}) &= \ph{-}0, \\
\phi^\prime(\overrightarrow{b}) &= \ph{-}Aj-Bi
\end{align*}
where $\phi^\prime$ denotes $\phi$ composed with projection onto the second coordinate of $H_1(H)$.

Recalling that performing surgery along $k$ turns $H$ into a solid torus $V$, we have a genus one Heegaard splitting of a lens space $L(p,q)$ with a core of $V$ represented by the generator $(0,1)$ of $H_1(H)$. In addition, $\phi(\overrightarrow{B})$ represents a meridian of the other solid torus $W$ of the genus one Heegaard splitting of $L(p,q)$ and $\phi(\overrightarrow{b})$ represents a longitude of $W$. It follows that $p= |\phi^\prime(\overrightarrow{B})| = |Aa+Bb|$.

Note that the multiple of $(0,1)$ which the core of $W$ represents in $H_1(L(p,q))$ determines $q$ mod $p$. So we want to evaluate $\phi^\prime(\overrightarrow{b})$ mod $p$. We will do this by observing that there is a relationship between $q$ and the homology class which the knot $k^\prime$ dual to $k$ represents in $L(p,q)$. Since $\phi(\overrightarrow{A})$ is homologous to $k^\prime$ in $L(p,q)$, we may take $\phi^\prime(\overrightarrow{A})$ as a representative of the homology class $k^\prime$ represents in $H_1(L(p,q))$.

Let $\lambda = \phi^\prime(\overrightarrow{A}) =-(Ar+Bs)$, and assume, for the moment, that neither $a$ nor $B$ is congruent to zero mod $p$. Then $a\lambda =-(Aar+Bas)$. Or, since, from the expression for $p$ above, $Aa\equiv -Bb$ mod $p$, $a\lambda \equiv -B(as-br)$ mod $p$. But, since $as-br=1$, $a\lambda =-B$ mod $p$. Then if $(a,p) = 1$, $\lambda \equiv -B/a$ mod $p$.

Next consider the equation $Ar+Bs=ai+bj$. Replacing $Ar+Bs$ with $-\lambda$, then multiplying both sides by $B$ and using the fact that $Bb \equiv -Aa$ mod $p$, yields $B\lambda \equiv a(Aj-Bi)$ mod $p$. Then using $\lambda \equiv -B/a$ mod $p$ and replacing $Aj-Bi$ with $\phi^\prime(\overrightarrow{b})$, we have $\phi^\prime(\overrightarrow{b}) \equiv -\lambda^2$ mod $p$ or equivalently, $a^2 \phi^\prime(\overrightarrow{b}) \equiv B^2$ mod $p$. Finally, the classification of lens spaces implies $\phi^\prime(\overrightarrow{b}) \equiv \pm q^{\pm1}$ mod $p$ . So $q \equiv \pm \lambda^{\pm2}$ mod $p$. Or $a^2 q \equiv \pm B^{\pm2}$ mod $p$. This is the first of the two claimed formulas for $q$. 

If $a$ or $B$ is congruent to zero mod $p$, then $b$ is not congruent to zero mod $p$ and $A$ is not congruent to zero mod $p$. So this derivation can be repeated, with appropriate changes, to obtain alternative expressions $\lambda \equiv A/b$ mod $p$ for $\lambda$ and $b^2 q \equiv \pm A^{\pm2}$ mod $p$ for q. Q.E.D.

\medskip
\centerline{\underline{Descriptions of the knots}.} 
\medskip

Those knots which lie in a solid torus and have non-trivial surgeries yielding a solid torus are a special case of the class of knots satisfying Theorem 1. The first six families of knots that we describe arise in this way. Next, we describe a family of knots each of which lives on a Seifert surface of a trefoil knot in $S^3$, and another family of knots each of which lives on a Seifert surface of a figure-eight knot in $S^3$. Then there are a final four families of knots satisfying Theorem 1, which we call ``\underline{sporadic}'', and we conclude by describing them.

\medskip
\centerline{\underline{The knots arising from knots in solid tori with surgeries yielding solid tori}.}
\medskip

I) \underline{Torus Knots}.

Let $k$ be the simple closed curve in $\partial H$ as illustrated by the R-R diagram in Fig. 2.

$A=1$, $(B,b)=1$, and $B \geq 2$.

\centerline{Insert Fig. 2 here.}

ii) \underline{Cables about Torus Knots}.

Let $k$ be the simple closed curve in $\partial H$ as illustrated by the R-R diagram in Fig. 3. 

$A=1$, $(B,b)=2$, and $B\geq 4$.

\centerline{Insert Fig. 3 here.}

III) 	Let $k$ be the simple closed curve in $\partial H$ obtained by twisting the simple closed curve $c$ $J$ times to the left about the curve $T$, as illustrated by the R-R diagram in Fig. 4.

\medskip
\centerline{Insert Fig. 4 here.}
\smallskip

One checks that, if $J \geq 0$, $n\geq 0$, $\varepsilon =\pm 1$ and $K$ is arbitrary, then $k$ satisfies Theorem 1. Let $A=J+1$, $B=(J+1)(2n+\varepsilon)$, $a=\pm 1$, and $b=-a(2\varepsilon(J+1)+BK)$. Then $[\overrightarrow{k}]$ has coordinates:
\[ [\overrightarrow{k}] = (A,B,a,b) \] in $H_1(\partial H)$. The constraints on $J$, $n$, and $\varepsilon$ become the following constraints on $A$, $B$, $a$, and $b$: $A>1$, there exists $\varepsilon =\pm 1$ such that $(B+\varepsilon)/A$ is an odd integer, $a=\pm 1$, and $b\equiv -2\varepsilon Aa$ mod $B$.

IV)	Let $k$ be the simple closed curve in $\partial H$ obtained by twisting the simple closed curve $c$ $J$ times to the right about the curve $T$, as illustrated by the R-R diagram in Fig. 5.

\medskip
\centerline{Insert Fig. 5 here.}

One checks that if $J \geq 0$, $n\geq 0$, $\varepsilon=\pm 1$, and $K$ is arbitrary, then $k$ satisfies Theorem 1. Let $A=2J+1$, $B=2n(J+1)+J\varepsilon)$, $a=\pm 1$, and $b=-a(\varepsilon(2J+1)+BK)$. Then $[\overrightarrow{k}]$ has coordinates:
$[\overrightarrow{k}] = (A,B,a,b)$  in $H_1(\partial H)$. The constraints on $J$, $n$, and $\varepsilon$ become the following constraints on $A$, $B$, $a$, and $b$:
$A>3$, there exists $\varepsilon=\pm 1$ such that $(2B+\varepsilon)/A$ is integral, $a=\pm 1$, and $b\equiv - \varepsilon Aa$ mod $B$.

V)	Let $k$ be the simple closed curve in $\partial H$ obtained by twisting the simple closed curve $c$ $J-1$ times to the right about the curve $T$, as illustrated by the R-R diagram in Fig. 6.

\medskip
\centerline{Insert Fig. 6 here.}

One checks that, if $J\geq 0$, $n \geq 0$, $\varepsilon=\pm 1$, and $K$ is arbitrary, then $k$ satisfies Theorem 1. Let $A=2J+1$, $B=n(2J+1)+J\varepsilon$, $a=\pm 1$, and $b=-a(\varepsilon(2J+1)+BK)$, then $[\overrightarrow{k}]$ has coordinates:
$[\overrightarrow{k}] = (A,B,a,b)$ in $H_1(\partial H)$. The constraints on $J$, $n$, and $\varepsilon$ become the following constraints of $A$, $B$, $a$, and $b$:
$A>1$ and odd, there exists $\varepsilon=\pm 1$ such that $(B-\varepsilon)/A$ is integral, $a=\pm 1$, and $b \equiv -\varepsilon Aa$ mod $B$. 

VI)	Let $\overrightarrow{k}$ be the curve obtained by twisting the curve c $J$ times to the left about the curve T as illustrated in Figure 7 below. 

\centerline{Insert Fig. 7 here.}

Then $[\protect \overrightarrow{k}] = (A,B,a,b)$, where $A>2$, $A$ is even, $B=2A+1$, and $b\equiv a(A-1)$ mod $B$ with $J \geq 1$.

\bigskip
\centerline{\underline{A family of knots obtained from a Seifert surface of a trefoil knot}.}

\bigskip
	Let $\overrightarrow{g}_1$ and $\overrightarrow{g}_2$ be simple closed curves embedded in $\partial H$ as in Figure 8.

\medskip
\centerline{Insert Fig. 8 here.}
\medskip

We mention that if $N$ is a regular neighborhood of $\overrightarrow{g_1}\cup\overrightarrow{g_2}$ in $\partial H$, then $\partial N$
is a trefoil knotted curve $c$ in $S^3$ and $N$ is a Seifert surface for $c$. In the chosen homology basis for $\partial H$, we have $[\overrightarrow{g_1}]=(1,0,-1,0)$ while $[\overrightarrow{g_2}]=(0,1,-1-1)$. Let $\overrightarrow{k}$ be a nonseparating
simple closed curve in $N$ with 
$[\overrightarrow{k}]=m[\overrightarrow{g_1}]+n[\overrightarrow{g_2}].$ 
Where, since $k$ is a simple closed curve, $(m,n)=1$. Note that the handlebody $H$ is homeomorphic to $N\times I$, and that the handlebody $H^\prime$ is also homeomorphic to $N\times I$. Thus any nonseparating simple closed curve in $N$ represents a free generator of both $\pi_1(H)$ and $\pi_1(H^\prime)$. It then follows from Theorem 1 that there is an integer surgery on $k$ yielding a lens space. Since $[\overrightarrow{k}]$ has coordinates: $[\overrightarrow{k}] = (m,n,-(m+n),-n)$
in $H_1(\partial H)$, it follows that the lens spaces obtained by integer surgeries on the knots in this family are of the form:
 $L(p,q) = L(m^2+mn+n^2,(m/n)^2).$

\underline{Remark}.
Observe that if $m$ or $n$ equals $0$, then $k$ is unknotted, while if $|m|$ or $|n|$ equals $1$, then $k$ is a torus knot. Hence, we may assume, without loss of generality, that $|m|$ and $|n|$ are greater than $1$.

\underline{Theorem 4}.
A lens space $L(p,q)$ is obtainable via an integer surgery on a knot embedded as a nonseparating simple closed curve in a genus one Seifert surface of a trefoil knot if and only if $p$ is odd, not divisible by $9$, every prime factor of $p$, not equal to $3$, is congruent to $1\mod 6$, and $q\equiv \pm((\pm-1\pm\sqrt{-3})/2)^{\pm 2}
\mod p$.

\underline{Proof}.
Consider the quadratic form $Ax^2+Bxy+Cy^2$. Recall from the theory of quadratic forms that the \underline{discriminant} of such a form is $B^2-4AC$ and that two quadratic forms in $x$ and $y$ are \underline{equivalent} if replacing $x$ with $ax+by$ and $y$ with $cx+dy$, where $ad-bc=1$, in one form carries that form into the other.

Note that $x^2+xy+y^2$ is a quadratic form with discriminant $-3$. If setting $x=a$ and $y=b$ yields $p$, with $(a,b)=1$, then we see that it is necessary that $-3$ be a quadratic residue $\mod 4p$. By applying the law of quadratic reciprocity, we see that the stated conditions on $p$ are exactly the conditions necessary for $-3$ to be a quadratic residue $\mod 4p$. On the other hand, it can be shown that these necessary conditions on $p$ are also sufficient conditions for $p$ to be expressible by some quadratic form of discriminant $-3$.
However, it is known that any two quadratic forms of discriminant $-3$ are equivalent [].
It follows from this fact that the conditions on $p$ are also sufficient to guarantee that $p$ can be expressed as $p=a^2+ab+b^2$ with $(a,b)=1$.

If $p=a^2+ab+b^2$ with $(a,b)=1$, let $\lambda \equiv a/b$ mod $p$. Then $\lambda$ satisfies the quadratic equation $\lambda^1+\lambda+1\equiv 0$ mod $p$. Thus $\lambda \equiv (-1\pm\sqrt{-3})/2$ mod $p$. Note that, for a given value of
$\sqrt{-3}$ mod $p$, the two possible values of $\lambda$ are reciprocals $\mod p$. Finally, by Lemma 3, $q \equiv\pm(\lambda)^{\pm 2}$ mod $p$. Q.E.D.

\bigskip
\centerline{\underline{A family of knots obtained from a Seifert surface of the figure-eight  knot}.}

\bigskip
Let $\overrightarrow{g_1}$ and $\overrightarrow{g_2}$ be simple closed curves embedded in $\partial H$ as illustrated by Fig. 9.

\medskip
\centerline{Insert Fig. 9 here.}

\medskip
Then, if $N$ is a regular neighborhood of $\overrightarrow{g_1}\cup\overrightarrow{g_1}$ in $\partial H$, $\partial N$ is a
figure-eight knotted curve $c$ in $S^3$ and $N$ is a Seifert surface for $c$. In the chosen homology basis for $\partial H$, we have $[\overrightarrow{g_1}]=(1,0,-1,0)$, while $[\overrightarrow{g_1}]=(0,1,-1,1)$. Let $\overrightarrow{k}$ be a nonseparating
simple closed curve in $N$ with
$ [\overrightarrow{k}]= m[\overrightarrow{g_1}]+n[\overrightarrow{g_2}]. $
Where, since $k$ is a simple closed curve, $(m,n)=1$. Note that the handlebody $H$ is homeomorphic to $N \boldsymbol{\times} I$, and the handlebody $H^\prime$ is also homeomorphic to $N \boldsymbol{\times} I$. Thus any nonseparating simple closed
curve in $N$ represents a free generator of both $\pi_1(H)$ and $\pi_1(H^\prime)$. It then follows from Theorem 1 that there is an integer surgery on $k$ yielding a lens space. Since $[\overrightarrow{k}]$ has coordinates:
$[\overrightarrow{k}]= (m,n,-(m+n),n)$ in $H_1(\partial H)$, it follows that the lens spaces obtained by integer surgeries on the knots in this family are of the form: $L(p,q)=L(|n^2-mn-m^2|,(m/n)^2).$

\underline{Remark}.
Observe that, as was the case for the family of knots based on a Seifert surface of a trefoil knot, if $m$ or $n$ equals $0$, then $k$ is unknotted, while if $|m|$ or $|n|$ equals $1$, then $k$ is a torus knot. Hence we may assume, without loss of generality, that $|m|$ and $|n|$ are greater than $1$.

The following theorem characterizes those lens spaces which can be obtained by surgeries on the knots in this family.

\underline{Theorem 5}.
A lens space $L(p,q)$ is obtainable via an integer surgery on a knot embedded as a nonseparating simple closed curve in a genus one Seifert surface of the figure-eight knot if and only if $p$ is odd, not divisible by $25$, every prime factor of $p$, not equal to $5$, is congruent to $\pm 1\mod 5$, and $q\equiv\pm((-1\pm\sqrt{5})/2)^{\pm 2}\mod p$.

\underline{Proof}.
The expression $b^2-ab-a^2$ corresponds to a quadratic form of discriminant $5$. Thus form can be analyzed in a fashion similar to that used for the form $a^2+ab+b^2$. We see that the stated conditions on $p$ are necessary in order that $p$ be representable by some quadratic form of discriminant $5$. However, it is not hard to check that any two quadratic forms of discriminant $5$ are equivalent. It follows from this fact that the conditions on $p$ are also sufficient to guarantee that $p$ can be represented as $p=b^2-ab-a^2$ with $(a,b)=1$.

If $p=b^2-ab-a^2$ with $(a,b)=1$, let $\lambda = a/b \mod p$. Then $\lambda$ satisfies the quadratic equation $\lambda^2+\lambda-1\equiv 0\mod p$. Thus $\lambda \equiv (-1\pm\sqrt{5})/2\mod p$. Note that, for a given value of $\sqrt{5}\mod p$, the two possible values of $\lambda$ are negative reciprocals mod $p$. Finally, by Lemma 3, $q=\pm(\lambda)^{\pm 2}\mod p$. Q.E.D.

\underline{Remark}.
Note that the construction of the proceeding two families works essentially because each knot in each family is a nonseparating simple closed curve embedded in a genus one Seifert surface of a fibered knot in $S^3$. If there were other genus one fibered knots besides the left- and right-handed trefoil and the figure-eight knot in $S^3$, or if the trefoil or the figure-eight knot had other inequivalent genus one Seifert surfaces, then we would expect to have other families of knots, with surgeries yielding lens spaces, obtained from curves embedded in these new Seifert surfaces. However, it is known, that the only fibered knots of genus one in $S^3$ are the trefoils and the figure-eight, and furthermore, that genus one Seifert surfaces of these knots are unique, up to isotopy []. So the families of knots obtained from Seifert surfaces of the trefoils and the figure-eight knot are the only such families that exist.

\bigskip
\centerline{\underline{The ``Sporadic'' Knots}}
\bigskip

Let $k$ be the simple closed curve in $\partial H$ obtained by twisting the simple closed curve $c$ $J$ times to the right about the curve $t$, as illustrated by the R-R diagram in Fig. 10.

\medskip
\centerline{Insert Fig. 10 here.}
\medskip

Note that in order for $k$ to be embedded in $\partial H$, the parameters $p$, $p'$, $m$, $m'$, $q$, $q'$, and $n$, $n'$ must be chosen so that $pm'-mp' = qn'-nq' = 1$.
Then observe that $\overrightarrow{k_\ast}$ represents $A^pB^n(A^mB^qA^mB^n)^J$ in $\pi_1(H)$, and $\overrightarrow{k_{\ast^ \prime}}$ represents  $a^{p^\prime}b^{n^\prime}(a^{m^\prime} b^{q^\prime}a^{m^\prime}b^{n^\prime})^J$ in $\pi_1(H^\prime)$. Now suppose the ordered set of parameters $(p,p^\prime,m,m^\prime, q,q^\prime,n,n^\prime)$ is replaced by one of the following four sets:
\begin{align*}
\text{a)} \quad &(\,1,\ph{-}1,\ph{-}2,\ph{-}3,\ph{-}1,-1,\ph{-}1,\ph{-}0\,) \\
\text{b)} \quad &(\,2,\ph{-}1,\ph{-}3,\ph{-}2,\ph{-}1,-1,\ph{-}1,\ph{-}0\,) \\
\text{c)} \quad &(\,4,-3,\ph{-}3,-2,\ph{-}1,\ph{-}0,\ph{-}1,\ph{-}1\,) \\
\text{d)} \quad &(\,3,-5,\ph{-}2,-3,\ph{-}1,\ph{-}0,\ph{-}1,\ph{-}1\,).
\end{align*}

Then, in each of these four cases, $\overrightarrow{k_\ast}$ and $\overrightarrow{k_{\ast^\prime}}$ represent free generators of $\pi_1(H)$ and
$\pi_1(H^\prime)$, respectively, provided $J\geq 0$. It then follows from Theorem 1 that, in each case, there is an integer surgery on $k$ yielding a lens space. 

The coordinates of $[\overrightarrow{k}]$ in $H_1(\partial H)$ corresponding to the preceding four possible sets of parameter values are:
\begin{align*}
\text{a)} \quad & (\,4J+1,\ph{1}2J+1,\ph{-}6J+1,-J\,) \\
\text{b)} \quad & (\,6J+2,\ph{1}2J+1,\ph{-}4J+1,-J\,) \\
\text{c)} \quad & (\,6J+4,\ph{1}2J+1,-4J-3,\ph{-}J+1\,) \\
\text{d)} \quad & (\,4J+3,\ph{1}2J+1,-6J-5,\ph{-}J+1\,).
\end{align*}

And so, by Lemma 3, the lens spaces obtained by integer surgeries on the knots in these four families are of the form:
\begin{align*}
\text{a)} \quad & L(\,22J^2+\ph{1}9J+\ph{1}1,\,(11J+2)^2\,) \\
\text{b)} \quad & L(\,22J^2+13J+\ph{1}2,\,(11J+3)^2\,) \\
\text{c)} \quad & L(\,22J^2+31J+11,\,(11J+8)^2\,) \\
\text{d)} \quad & L(\,22J^2+35J+14,\,(11J+9)^2\,).
\end{align*}

This completes the list of knots satisfying Theorem 1 known to us. We have not been able to prove that this list is complete. However, there is some computational evidence that it is complete. In particular, John Rosenberg kindly wrote a computer program for us which examined one-bridge knots, in lens spaces with fundamental groups of order less than 1000, for integer surgeries yielding $S^3$, without finding any knots not on our list.

A presumably more interesting and more difficult problem than characterizing double-primitive knots on genus two Heegaard surfaces in $S^3$ is that of characterizing all knots in $S^3$ with Dehn-surgeries yielding lens spaces. Then there is the following natural question.

\underline{Question}.
If $k$ is a knot in a lens space $M$ such that Dehn-surgery on $k$ yields $S^3$, must $k$ be a\, 0 or 1-bridge knot in $M$?

The following table includes all lens spaces $L(p,q)$ with $p\leq 500$ which are obtainable via surgeries on knots in our list,
except, to keep the size of the table managable,  those lens spaces which are obtainable only by surgeries on torus knots and cables about
torus knots are not shown. Included with each lens space $L(p,q)$, in column $\lambda$, is a list of the homology classes
in $L(p,q)$ that are represented by 1-bridge knots which have Dehn-surgeries that yield $S^3$. For example, the table shows that $L(37,10)$ contains two such 1-bridge knots.

In order to check whether a lens space $L(p,q)$ with $p>500$ is obtainable from a knot on these lists, some computation is necessary. Checking whether $L(p,q)$ is obtainable from a torus knot, a cable about a torus knot, a knot
based on a Seifert surface for a trefoil or figure-eight knot or a sporadic knot, is straightforward. To determine whether $L(p,q)$ is obtainable from one of the remaining knots, one can do the following: The relationship $p=Aa+Bb$ holds with
$0<2A\leq B$ and $|a|=1$. Thus the value of $B$ determines the value of $A$, $a$, and $b$. It then suffices to check those values of $B$ such that $q\equiv \pm(B)^{\pm 2}\mod p$, for $1< B \leq\max\{8,2p/5\}$, and to check whether the set of values $\{A,B,a,b\}$ obtained satisfies the constraints for any of the knots in our list.

\medskip
\centerline{Insert Table of Lens Spaces here.}
\medskip

\underline{References}.

\newpage
\setlength{\tabcolsep}{5.30pt}
\renewcommand{\arraystretch}{0.95}

\footnotesize

\begin{center}
\begin{tabular}{||c|c|c||c|c|c||c|c|c||c|c|c||c|c|c||c|c|c||c|c|c||}\hline
$p$&$q$ &$\lambda$&$p$&$q$ &$\lambda$&$p$&$q$ &$\lambda$&$p$&$q$ &$\lambda$&$p$&$q$ &$\lambda$&$p$&$q$ &$\lambda$&$p$&$q$ &$\lambda$ \\ \hline
18&5&5&&&25&142&53&35&189&41&22&226&49&7&269&61&21&&135&135 \\ 
19&7&7&94&35&23&143&25&5&&62&62&&69&35&&70&71&302&83&41 \\ 
27&8&8&95&39&18&144&47&47&&67&32&&85&57&270&59&31&&109&13 \\ 
30&11&7&&41&42&145&51&52&190&71&47&227&49&7&& 89&89&&113&75 \\ 
31&11&12&97&35&35&&59&28&191&56&18&&84&10&& 101&67&305&121&62 \\  
& 12&7&98&27&13&146&55&37&&87&88&229&64&8&& 109&53&306&101&101 \\ 
32&7&5&&37&25&147&41&20&192&71&11&&80&81&271&28&28&&115&77 \\ 
34&13&9&99&29&13&&62&8&193&25&5&&94&18&& 75&14&307&25&5 \\ 
36&11&11&&&16&&67&67&&81&9&&&94&273&43&20&&54&19 \\ 
37&10&8&&32&32&148&41&19&&84&84&230&49&13&& 100&10&308&87&45 \\  
&&10&100&39&19&149&39&40&194&73&49&&91&47&& &100&309&46&46 \\ 
39&16&16&101&21&22&150&59&29&195&79&38&231&67&32&274&81&9&311&57&58 \\ 
43&12&5&&30&8&151&26&27&196&55&27&232&25&5&& 103&69&&115&14 \\ 
45&14&14&103&46&46&&32&32&&75&9&233&89&12&275&49&7&312&49&7 \\  
& 19&8&105&41&22&&45&14&&&11&234&35&19&& 64&8&313&49&7 \\ 
46&17&11&107&25&5&153&50&50&198&37&17&&77&77&& 109&54&&71&23 \\ 
47&13&9&&41&24&&55&26&&65&65&235&51&24&276&49&7&&98&98 \\ 
49&18&18&108&23&13&155&46&8&199&37&19&&54&17&& 73&35&&119&11 \\ 
50&19&9&&35&35&&61&32&&55&12&236&55&29&277&59&13&315&41&22 \\  
& &13&109&45&8&157&25&5&&60&61&237&55&55&& 81&9&&101&52 \\ 
54&17&17&&&45&&36&11&&92&92&&64&8&& &14&&104&104 \\ 
55&21&12&110&41&27&&46&16&200&79&39&&&10&& 116&116&317&85&42 \\ 
57&16&5&112&31&9&158&37&11&201&37&37&238&89&59&278&85&11&&121&11 \\ 
59&24&25&&&17&&59&39&205&36&13&239&70&13&279&83&11&&&14 \\ 
61&13&13&114&43&29&161&45&24&&46&47&&71&36&& 89&46&318&25&5 \\  
& 16&17&115&34&9&&61&30&&61&11&240&71&13&& 92&92&&119&79 \\ 
62&23&15&&&14&162&35&17&&81&42&241&45&14&280&81&39&319&59&28 \\ 
63&17&10&116&45&7&&53&53&206&63&11&&50&51&& 111&57&&138&139 \\  
& 20&20&&&25&&61&41&&77&51&242&43&21&281&36&37&&149&150 \\ 
64&23&19&117&38&38&163&44&10&207&25&5&&91&61&282&25&5&320&129&63 \\ 
66&25&17&&43&7&&58&58&&65&34&243&25&5&283&44&44&322&121&81 \\ 
67&18&7&&&20&&62&15&&68&68&&53&26&& 64&8&323&89&15 \\  
& 29&29&118&25&5&165&29&16&209&79&80&&77&40&& 78&19&324&49&7 \\ 
68&19&5&119&50&13&&49&7&&80&9&&80&80&& 104&20&&71&35 \\  
& &7&120&49&23&&&8&&&17&245&69&34&286&51&27&&107&107 \\ 
70&29&13&121&35&36&166&49&7&210&59&31&&99&48&& 107&71&325&49&7 \\ 
71&21&11&&46&14&167&46&11&&79&53&247&68&68&287&45&23&&129&64 \\  
& 26&20&125&49&24&168&25&5&211&31&32&&87&87&&62&15&326&97&11 \\ 
72&23&23&126&41&41&169&22&22&&45&16&250&99&49&288&95&95&327&100&10 \\ 
73&27&10&&47&31&170&69&33&&64&8&251&116&117&& 119&11&&&11 \\ 
75&29&11&127&19&19&171&50&11&212&57&27&252&83&83&& &13&&154&154 \\  
& &14&&27&10&&53&28&&81&9&253&74&20&290&81&9&329&95&46 \\ 
78&29&19&&29&15&&56&56&213&59&19&254&95&63&& 109&73&330&131&67 \\ 
79&23&23&128&47&7&172&39&21&214&49&7&255&101&52&291&61&61&331&31&31 \\  
& 28&29&&&9&&63&31&215&49&7&257&25&5&292&111&11&&75&16 \\  
& 29&7&129&49&7&173&64&8&216&49&23&&59&14&293&25&5&&76&14 \\ 
80&31&7&&&49&&&10&&71&71&258&97&65&& 64&8&&115&116 \\  
& &17&130&49&33&174&65&43&217&25&25&259&73&38&& 81&9&332&25&5 \\ 
81&26&26&&51&27&175&69&34&&67&67&&100&100&294&83&41&333&110&110 \\  
& 31&14&131&50&9&177&49&7&218&25&5&&121&121&295&64&22&&115&56 \\ 
82&23&5&132&25&5&&&17&219&50&13&261&86&86&& 91&92&334&125&83 \\  
& 31&21&133&30&30&178&49&7&&64&8&&91&44&& 119&58&335&71&34 \\ 
83&19&8&&36&13&&67&45&&&64&&100&19&297&67&32&337&63&20 \\ 
84&25&11&&39&18&179&73&74&220&89&43&263&49&7&& 98&98&&91&10 \\ 
85&26&11&135&31&14&180&59&59&222&83&55&&71&10&& 103&50&&94&22 \\ 
89&34&12&&41&22&&71&37&223&39&39&264&49&7&300&71&37&&128&128 \\ 
90&29&29&&44&44&181&48&48&&66&17&265&49&26&& 119&59&338&51&25 \\ 
91&16&16&137&37&10&&70&16&225&74&74&267&79&13&301&45&16&&127&85 \\  
& 27&8&139&42&42&182&25&5&&79&38&268&25&5&& 64&8&339&64&8 \\ 
93&25&5&&62&63&&53&25&&89&44&&99&13&& 79&79&340&89&43 \\ \hline
\end{tabular}
\end{center}

\newpage
\setlength{\tabcolsep}{4.40pt}
\renewcommand{\arraystretch}{0.95}

\begin{center}
\begin{tabular}{||c|c|c||c|c|c||c|c|c||c|c|c||c|c|c||c|c|c||c|c|c||}\hline
$p$&$q$ &$\lambda$&$p$&$q$ &$\lambda$&$p$&$q$ &$\lambda$&$p$&$q$ &$\lambda$&$p$&$q$ &$\lambda$&$p$&$q$ &$\lambda$&$p$&$q$ &$\lambda$ \\ \hline
341&74&24&362&49&7&&150&151&411&49&7&&163&109&&118&118&475&64&8 \\   
&79&80&363&65&32&391&53&27&&64&8&435&91&44&455&81&9&&189&94 \\  
&100&21&&98&10&392&111&55&&76&26&&94&23&&129&66&476&137&67 \\  
342&61&29&364&87&45&393&25&5&413&121&11&436&81&9&&181&92&477&158&158 \\   
&113&113&365&64&8&395&159&78&&&16&&129&13&457&25&5&&163&80 \\  
343&25&5&&69&36&&186&187&414&137&137&437&100&10&&133&133&478&179&119 \\   
&54&17&366&137&91&396&131&131&&155&103&&169&16&&169&13&479&227&228 \\   
&97&48&367&83&83&397&34&34&415&54&19&&&13&&72&23&480&191&97 \\   
&131&12&368&25&5&398&149&99&&116&24&&136&14&&133&18&481&100&100 \\  
345&139&68&369&80&17&399&85&13&416&95&17&438&121&11&459&49&7&&211&211 \\  
347&61&16&&122&122&&121&11&417&181&181&439&68&69&&55&26&482&25&5 \\   
&64&8&&127&62&&&121&418&25&5&&78&19&&103&50&&181&121 \\   
&108&14&370&139&93&&163&163&&157&105&&81&9&&149&76&484&87&43 \\  
349&106&22&&149&73&400&119&11&419&95&18&&171&171&&152&152&485&64&8 \\   
&122&122&371&81&9&&159&79&420&169&83&440&81&39&460&49&7&486&85&41 \\   
&142&143&373&49&7&401&110&111&421&63&22&441&125&62&461&100&19&&107&53 \\ 
350&131&87&&88&88&&111&17&&64&8&&146&146&462&173&115&&161&161 \\   
&139&69&&100&10&402&151&101&&74&16&&151&74&463&63&20&487&91&17 \\  
351&77&40&374&49&7&403&61&32&&109&110&442&69&33&&100&10&&105&52 \\   
&113&58&&81&9&&64&8&422&49&7&443&25&5&465&86&28&&144&12 \\   
&116&116&&135&13&&87&87&423&49&7&&79&27&&89&46&&232&232 \\  
352&97&15&375&149&74&&191&191&&137&70&&129&25&466&129&17&488&73&19 \\  
354&133&89&377&144&12&404&105&51&&140&140&445&179&88&&175&117&489&58&58 \\  
355&61&62&378&59&31&405&89&44&425&169&84&&186&187&467&64&8&490&139&69 \\   
&66&17&&85&41&&134&134&427&74&74&446&95&13&&166&16&491&72&73 \\   
&81&9&&109&53&&139&68&&100&10&&167&111&&193&14&492&119&61 \\   
&99&16&&125&125&&161&82&&123&60&447&121&11&468&25&5&493&25&5 \\   
&141&72&379&51&51&406&93&47&&135&135&448&121&11&&121&59&&64&8 \\  
357&25&5&&105&22&&115&59&428 &103&53&449&164&165&&155&155&494&185&123 \\   
&64&8&380&71&17&407&25&5&429&64&8&450&59&29&469&37&37&495&161&82 \\   
&101&52&&151&77&&73&38&430&161&107&&149&149&&71&34&&164&164 \\ 
358&81&9&382&25&5&408&47&25&&171&87&&169&113&&163&163&&199&98 \\  
359&82&21&&143&95&&169&13&431&89&90&&179&89&470&189&93&498&187&125 \\   
&104&105&386&145&97&409&53&53&432&25&5&451&46&47&471&49&7&&193&17 \\  
360&49&23&387&125&64&&120&23&&95&49&&103&17&&169&13&499&89&29 \\   
&119&119&&128&128&&121&11&&143&143&&125&24&&&169&&139&139 \\  
361&41&42&388&163&13&&&13&433&98&26&&156&157&472&49&7&&223&224 \\   
&49&7&&&15&&128&129&&198&198&452&81&9&473&100&10&&&\\   
&68&68&389&91&46&410&49&7&434&121&11&453&104&14&&181&14&&&\\ \hline
\end{tabular}
\end{center}

\begin{figure}[h]
\centering
\includegraphics[width = .50\textwidth]{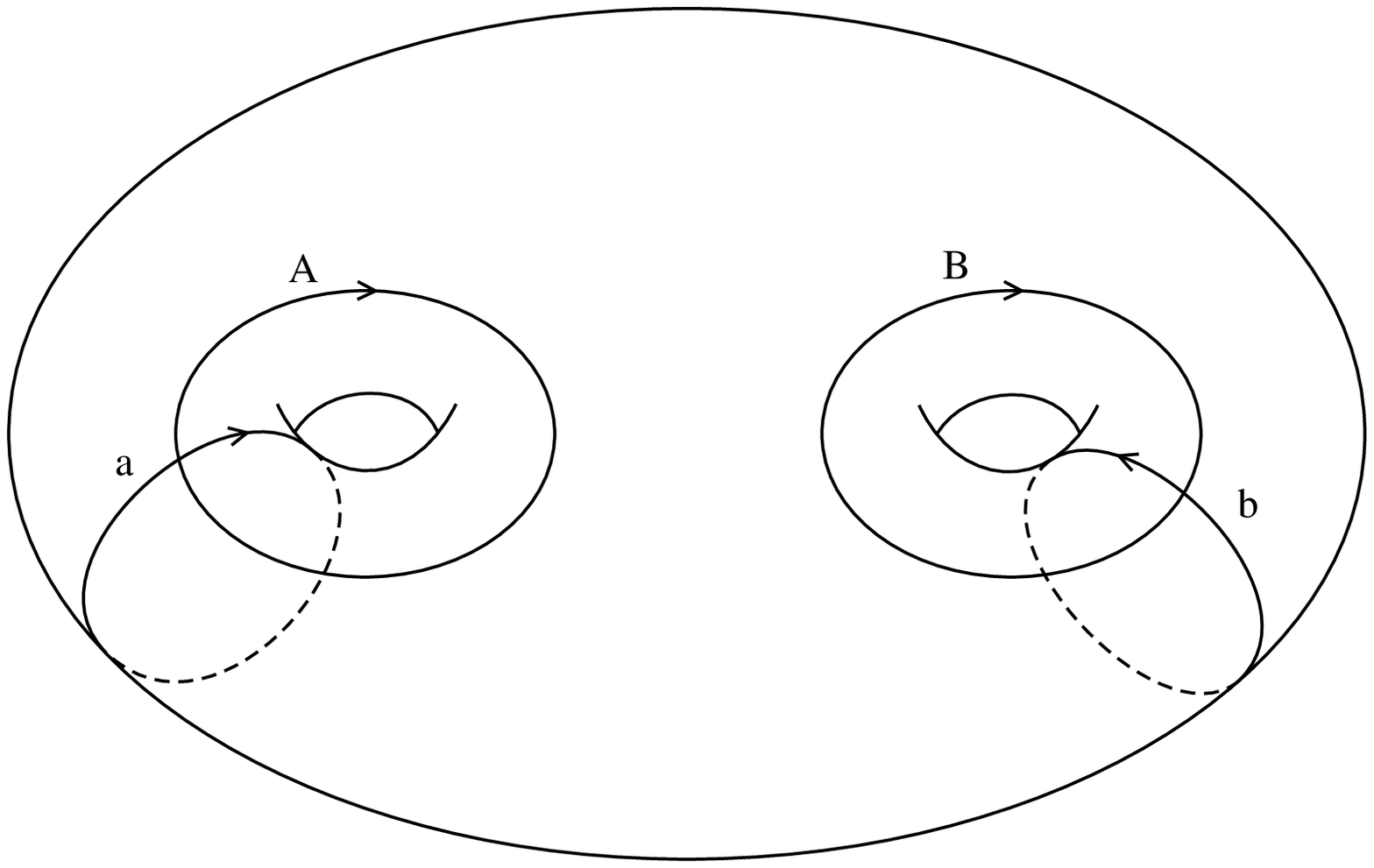}
\caption{A genus two handlebody $H$ standardly embedded in $S^3$ and a standard basis for $H_1(\partial H)$ of the Heegaard surface $\partial H$, with $\protect \overrightarrow{a}$, $\protect \overrightarrow{b}$ bounding disks in $H$ and $\protect \overrightarrow{A}$, $\protect \overrightarrow{B}$ bounding disks in the complementary genus two handlebody $H^\prime$.}
\label{Fi:1}
\end{figure}

\begin{figure}[h]
\centering
\includegraphics[width = .50\textwidth]{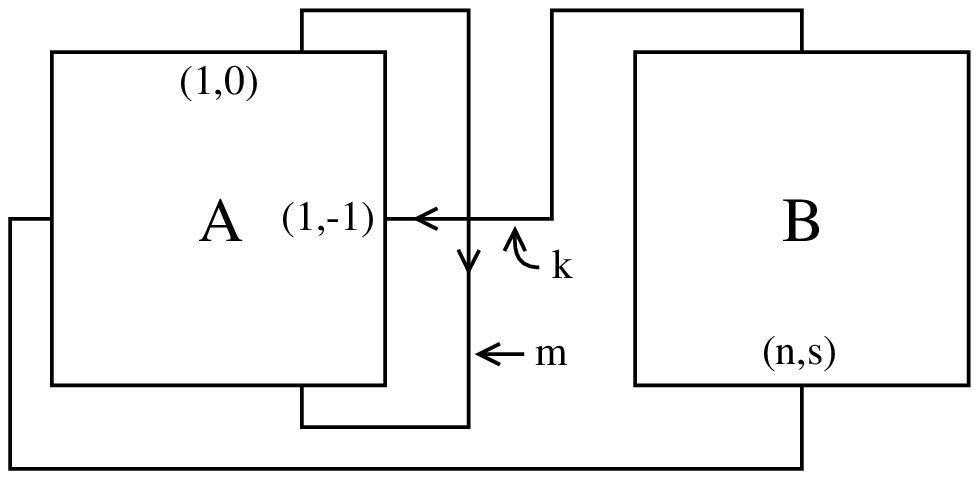}
\caption{
\centerline{\quad}
\centerline{\textbf{Diagrams of Type I: Torus Knots.}}
Here $\protect \overrightarrow{k}_*$ = $AB^n$ and $[\protect \overrightarrow{k}]$ = $(A,B,a,b)$, where $A=1$, $B = n$, $a = \pm 1$ and $b = s$, with $(n,s) = 1$, and $n \geq 2$; so $(B,b)=1$, and $B \geq 2$.}
\label{Fi:2}
\end{figure}

\begin{figure}[h]
\centering
\includegraphics[width = .50\textwidth]{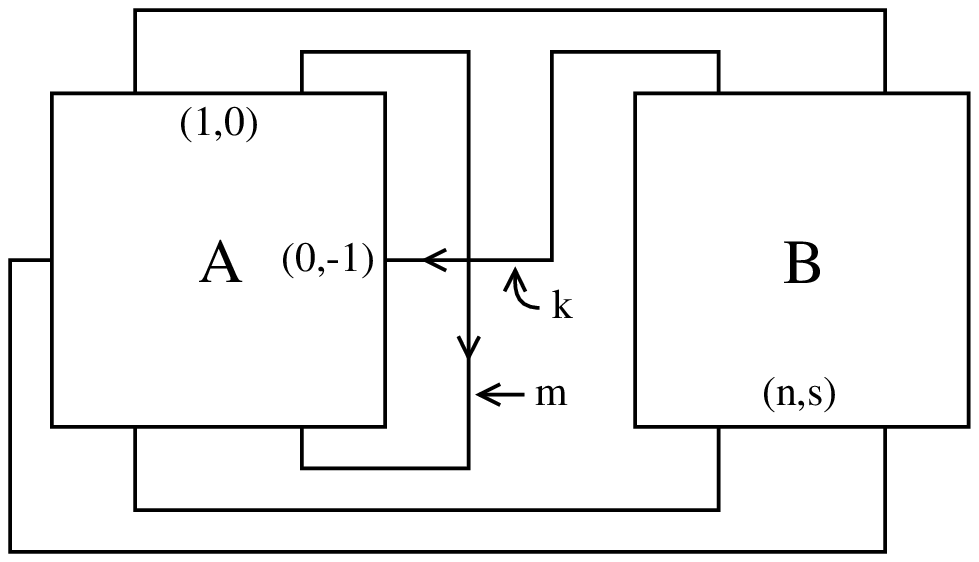}
\caption{
\centerline{\quad}
\centerline{\textbf{Diagrams of Type II: Cables of Torus Knots.}} 
Here $\protect \overrightarrow{k}_* = B^nAB^n$ and $[\protect \overrightarrow{k}] = (A,B,a,b)$, where $A=1$, $B = 2n$ and $b = 2s$, with $(n,s) = 1$, and $n \geq 2$; so $(B,b)=2$, and $B\geq 4$.}
\label{Fi:3}
\end{figure}

\begin{figure}[h]
\centering
\includegraphics[width = .60\textwidth]{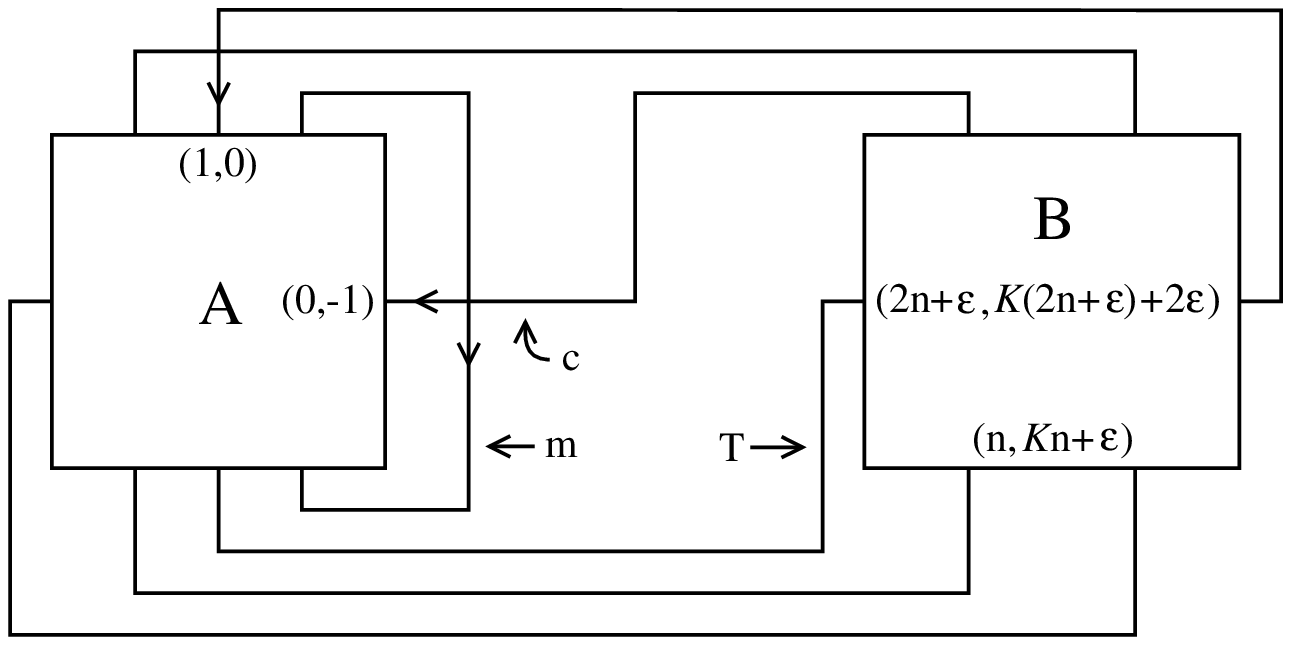}
\caption{
\centerline{\quad}
\centerline{\textbf{Diagrams of Type III Knots.}} 
Here $\protect \overrightarrow{k}$ is the curve obtained by twisting the curve c $J$ times to the left about the curve T. Then, $\protect \overrightarrow{k}_* = B^n(AB^{2n+\varepsilon})^JAB^n$ and $[\protect \overrightarrow{k}] = (A,B,a,b)$, where $A=J+1$, $B = (2n+\varepsilon)A - \varepsilon$, $a=\pm 1$, and $b \equiv -2a\varepsilon A$ mod $B$ with $J \geq 1$, $K$ arbitrary, $n \geq 1$, $\varepsilon = \pm 1$, $\varepsilon = 1$ if $n = 1$.}
\label{Fi:4}
\end{figure}

\begin{figure}[h]
\centering
\includegraphics[width = .60\textwidth]{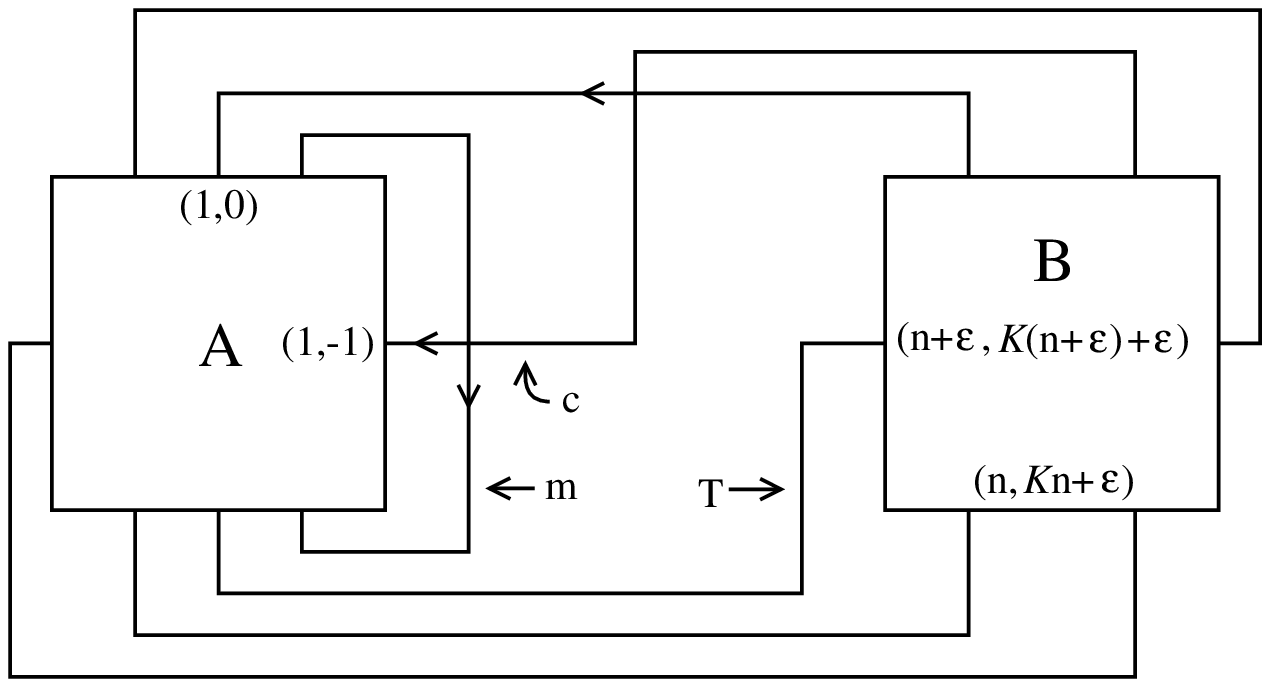}
\caption{
\centerline{\quad}
\centerline{\textbf{Diagrams of Type IV Knots.}} 
Here $\protect \overrightarrow{k}$ is the curve obtained by twisting the curve c $J$ times to the right about the curve T. Then, $\protect \overrightarrow{k}_* = AB^n(AB^{n+\varepsilon}AB^{n})^J$ and $[\protect \overrightarrow{k}] = (A,B,a,b)$, where $A=2J+1$, $B=nA+J\varepsilon$, $a=\pm 1$, and $b \equiv -a\varepsilon A$ mod $B$ with $J \geq 2$, $K$ arbitrary, $n \geq 2$, $\varepsilon = \pm 1$, $\varepsilon = 1$ if $n = 2$.} \label{Fi:5}
\end{figure}

\begin{figure}[h]
\centering
\includegraphics[width = .60\textwidth]{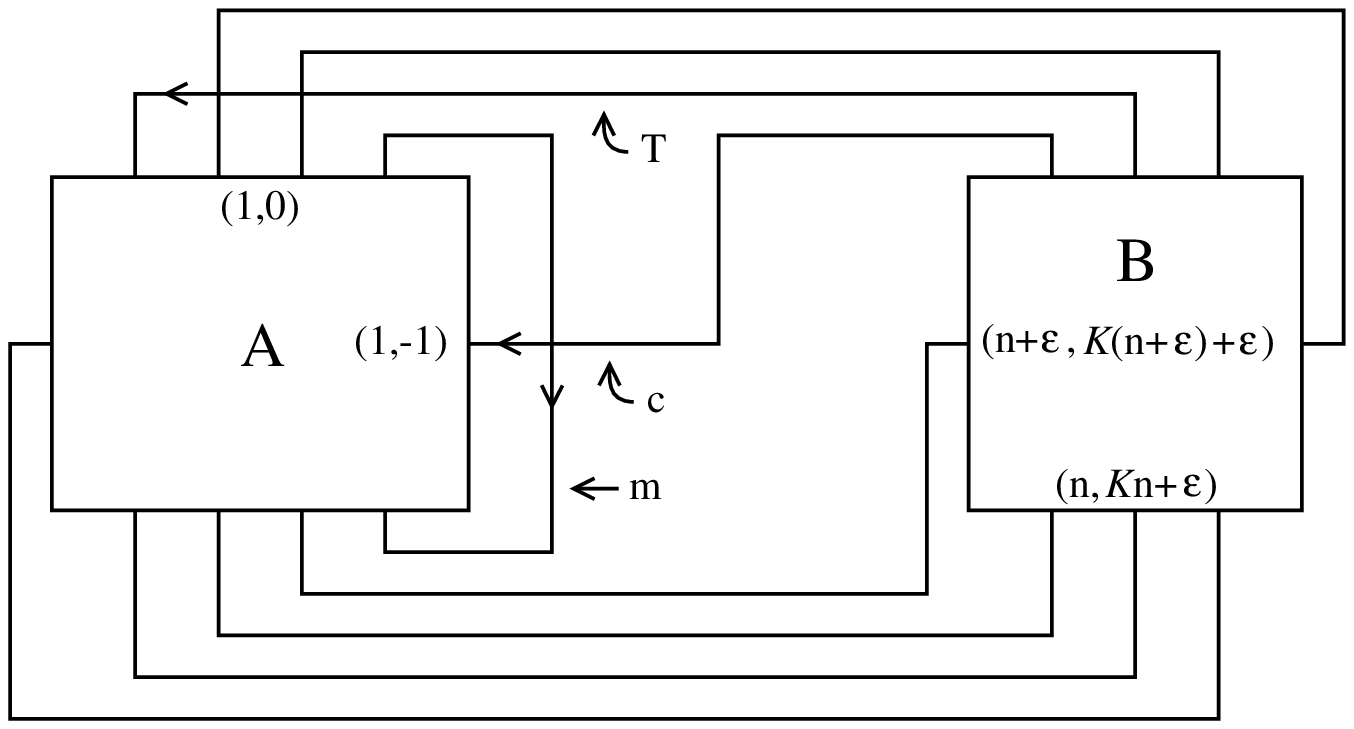}
\caption{
\centerline{\quad}
\centerline{\textbf{Diagrams of Type V Knots.}} Here $\protect \overrightarrow{k}$ is the curve obtained by twisting the curve c $J$ times to the right about the curve T. Then, $\protect \overrightarrow{k}_* = (AB^n)^{J+1}AB^{n+\varepsilon}(AB^{n})^{J+1}$ and $[\protect \overrightarrow{k}] = (A,B,a,b)$, where $A=2J+3$, $B=nA+\varepsilon$, $a=\pm 1$, and $b \equiv -a\varepsilon A$ mod $B$ with $J \geq 0$, $K$ arbitrary, $n \geq 2$, $\varepsilon = \pm 1$, $\varepsilon = 1$ if $n = 2$.}
\label{Fi:6}
\end{figure}

\begin{figure}[h]
\centering
\includegraphics[width = .60\textwidth]{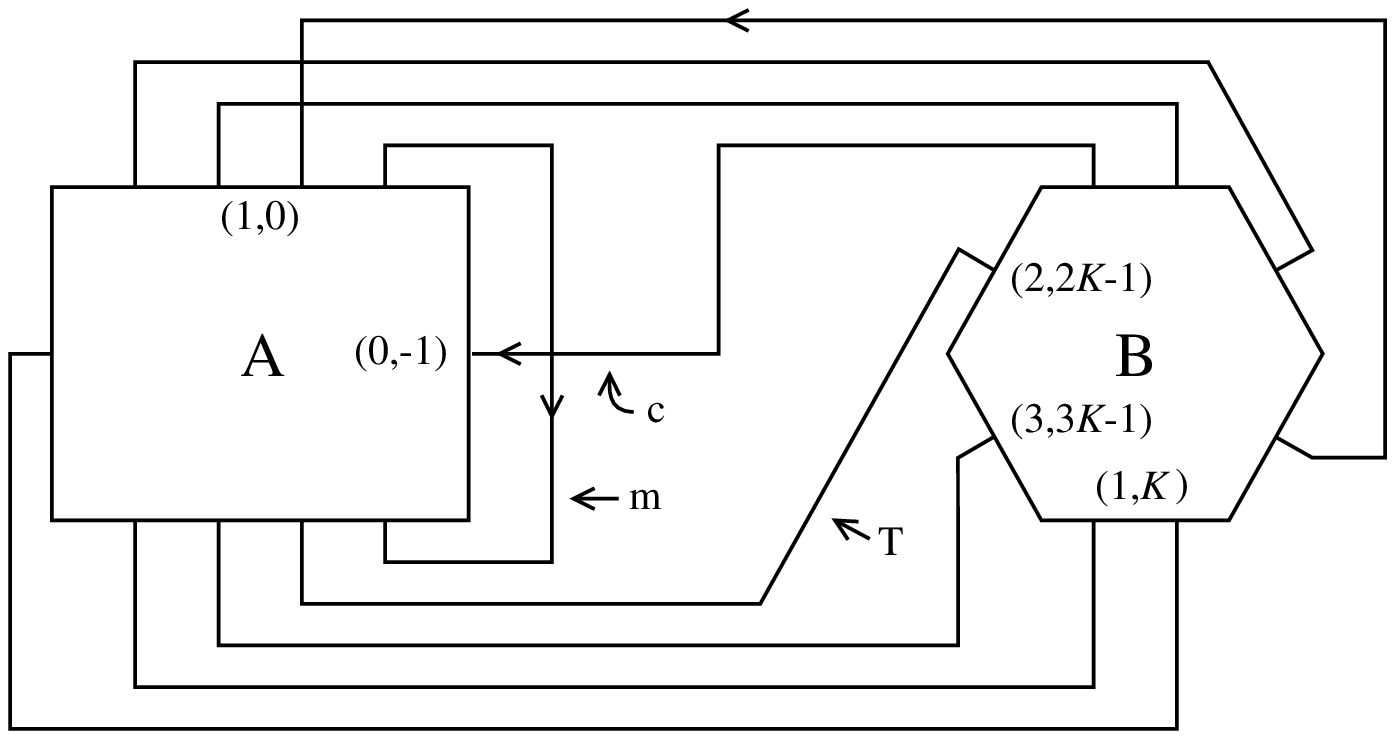}
\caption{
\centerline{\quad}
\centerline{\textbf{Diagrams of Type VI Knots.}} Here $\protect \overrightarrow{k}$ is the curve obtained by twisting the curve c $J$ times to the left about the curve T. Then, $\protect \overrightarrow{k}_* = B(AB^2)^JAB^3(AB^2)^JAB$ and $[\protect \overrightarrow{k}] = (A,B,a,b)$, where $A>2$, $A$ is even, $B=2A+1$, and $b\equiv a(A-1)$ mod $B$ with $J \geq 1$ and $K$ arbitrary.}
\label{Fi:7}
\end{figure}

\begin{figure}[h]
\centering
\includegraphics[width = .55\textwidth]{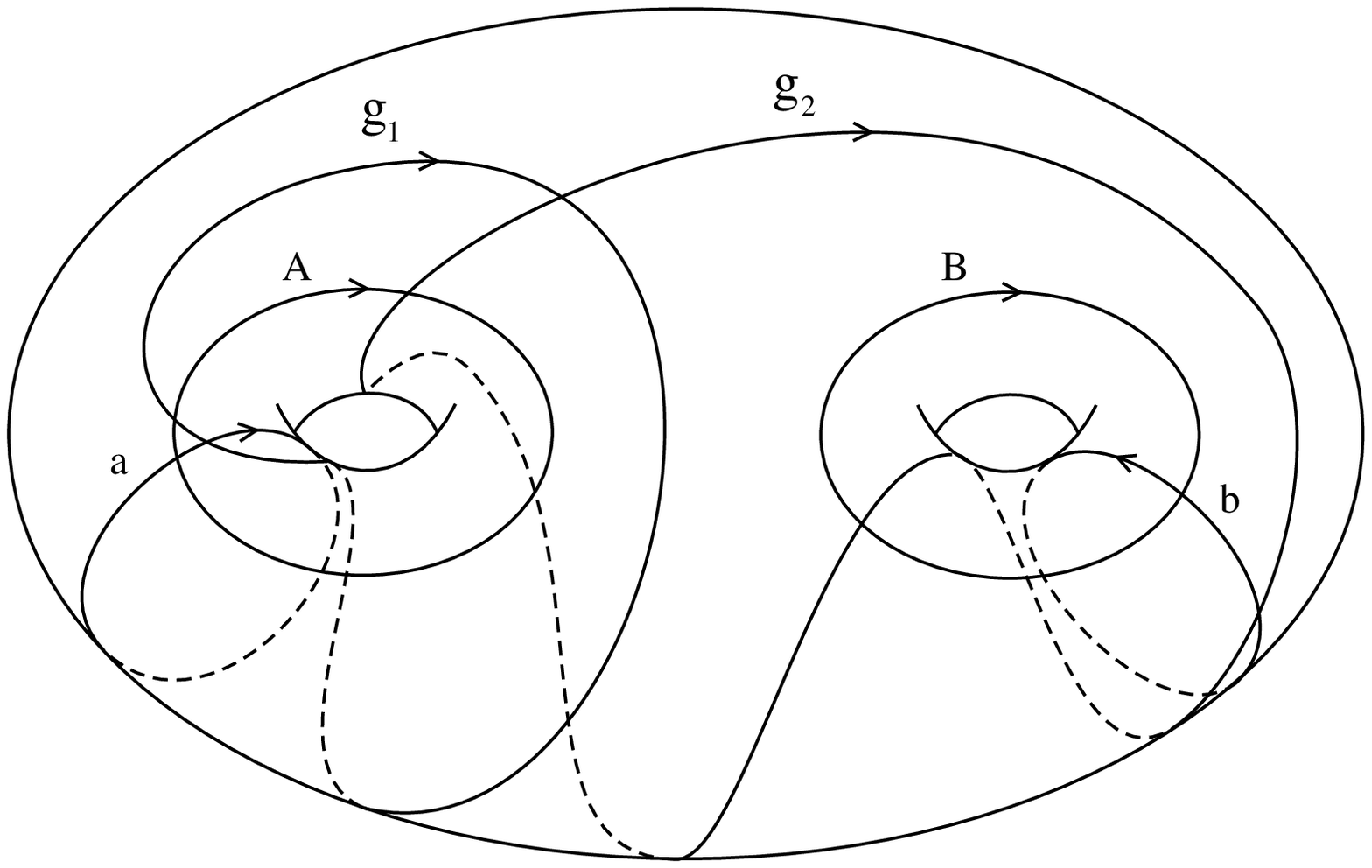}
\caption{Two simple closed curves $\protect \overrightarrow{g}_1, \protect \overrightarrow{g}_2$ placed in the boundary of a standardly embedded handlebody $H$ in $S^3$ so that a regular neighborhood of the pair $\protect \overrightarrow{g}_1, \protect \overrightarrow{g}_2$ in $\partial H$ is a fiber of a trefoil knot.}
\label{Fi:8}
\end{figure}

\begin{figure}[h]
\centering
\includegraphics[width = .55\textwidth]{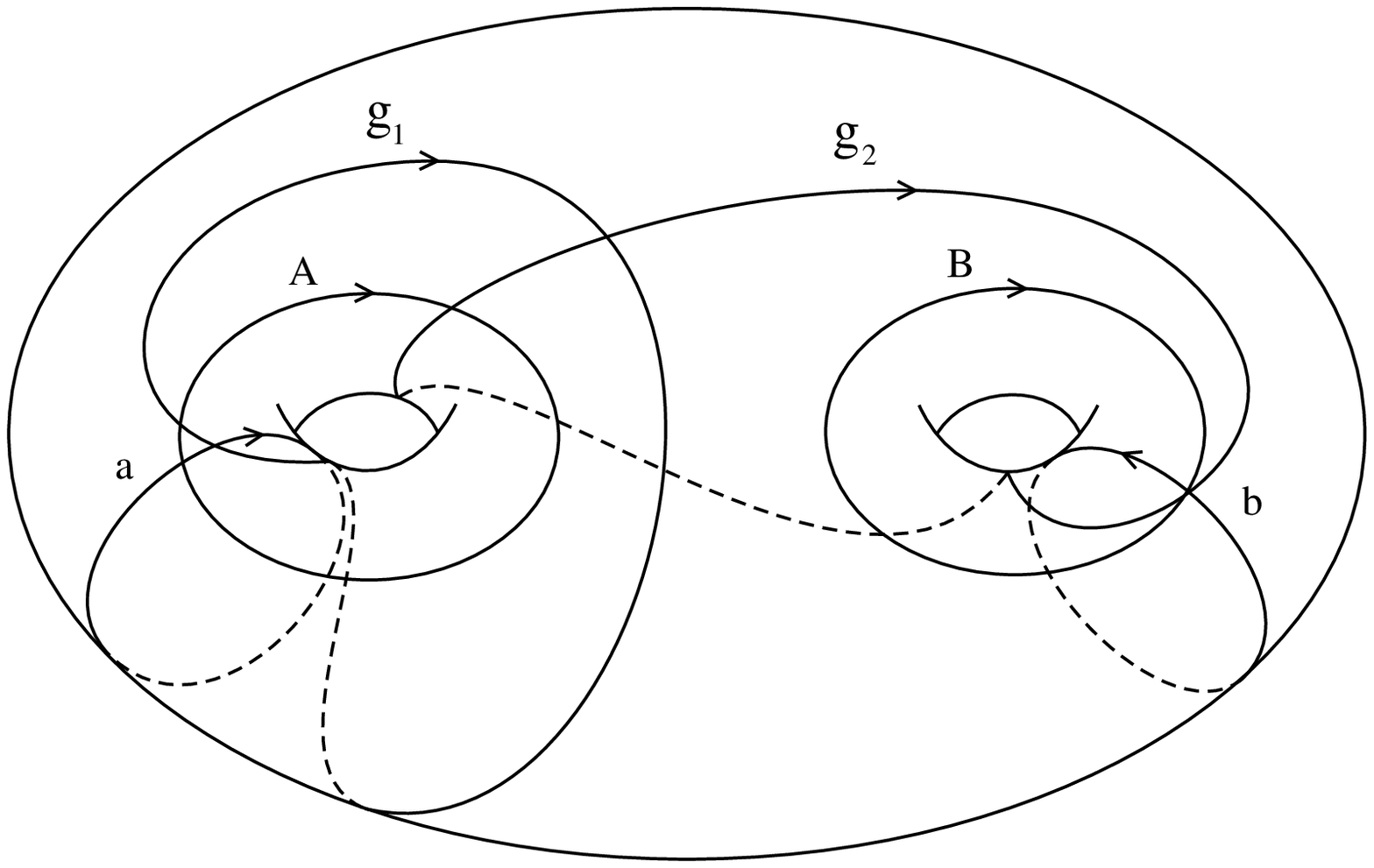}
\caption{Two simple closed curves $\protect \overrightarrow{g}_1, \protect \overrightarrow{g}_2$ placed in the boundary of a standardly embedded handlebody $H$ in $S^3$ so that a regular neighborhood of the pair $\protect \overrightarrow{g}_1, \protect \overrightarrow{g}_2$ in $\partial H$ is the fiber of the figure-eight knot.}
\label{Fi:9}
\end{figure}

\begin{figure}[h]
\centering
\includegraphics[width = .50\textwidth]{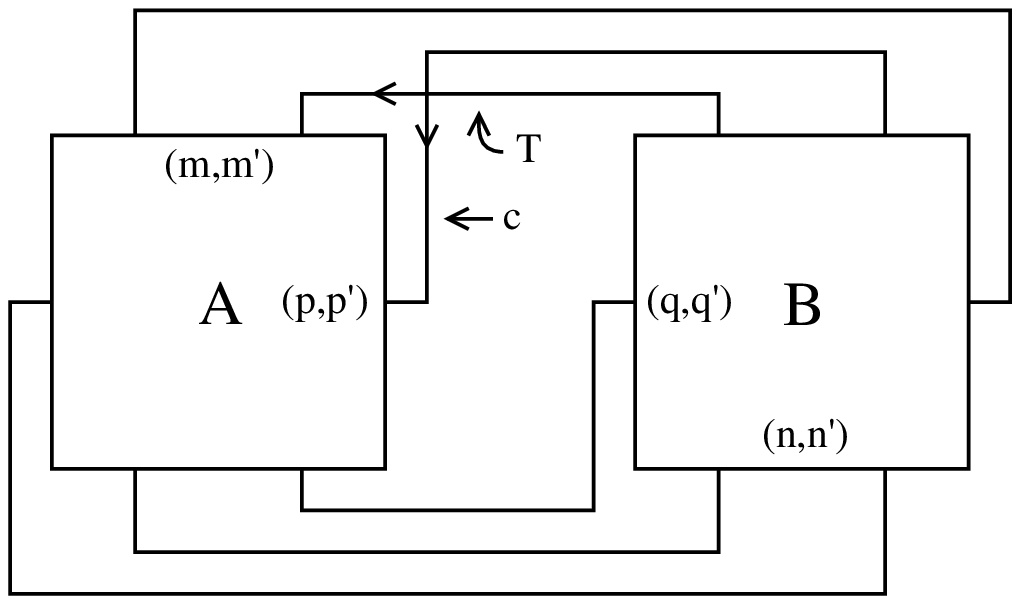}
\caption{
\centerline{\quad}
\centerline{\textbf{The ``Sporadic'' Knots.}} 
Suppose $\protect \overrightarrow{k}$ is the curve obtained by twisting the curve $c$ in the figure $J$ times to the right about the curve T. Then  $\protect \overrightarrow{k_\ast}$ represents $A^pB^n(A^mB^qA^mB^n)^J$ in $\pi_1(H)$, while $\protect \overrightarrow{k_{\ast^ \prime}}$ represents  $a^{p^\prime}b^{n^\prime}(a^{m^\prime} b^{q^\prime}a^{m^\prime}b^{n^\prime})^J$ in $\pi_1(H^\prime)$. If the 8-tuple of exponents $(p,p^\prime, \dots ,n,n^\prime)$ is replaced by $(\,1,1,2,3,1,-1,1,0\,)$, $(\,2,1,3,2,1,-1,1,0\,)$, $(\,4,-3,3,-2,1,0,1,1\,)$, or $(\,3,-5,2,-3,1,0,1,1\,)$, and $J \geq 0$, then $\protect \overrightarrow{k}$ becomes a double-primitive.}
\label{Fi:10}
\end{figure}

email address: jberge@charter.net
\end{document}